\documentstyle[12pt]{amsart}
 
\newcommand{\Nnn}{{\Bbb N}} 
\newcommand{\Cnn}{{\Bbb C}}

\newcommand{\ler}{left-to-right }
\newtheorem{theorem}{Theorem}

\newtheorem{lemma}{Lemma}
\newtheorem{proposition}{Proposition}
\newtheorem{definition}{Definition}
\newtheorem{corollary}{Corollary}
\newtheorem{example}{Example}
\newcommand{\deb}{D_n^{\beta(0,1)}}
\newcommand{\btree}{$\beta(0,1)$-tree }
\newenvironment{proof}{\noindent {\bf Proof:}}{{\qed}}
 \newcommand{\btrees}{$\beta(0,1)$-trees }

\newcommand{\vanish}[1]{}
\newcommand{\prek}{$P$-recursive }
 
\begin{document}
 
\title[Exact enumeration of 1342-avoiding permutations]
{Exact enumeration of 1342-avoiding permutations\\
A close link with labeled trees and planar maps}
 
\author{Mikl\'os B\'ona}
\address{\hskip-\parindent
        Department of Mathematics \\
        Massachusetts Institute of Technology \\
        Cambridge, MA 02139}
\email{bona@@math.mit.edu}
\thanks{Research at MSRI is supported in part by NSF grant DMS-9022140}
 
 \begin{abstract} Solving the first nonmonotonic, longer-than-three
instance of a 
classic enumeration problem,  we obtain the generating  function $H(x)$
of all 1342-avoiding permutations of length $n$ as well as an {\em exact}
formula for their number $S_n(1342)$.   While
achieving this, we bijectively prove that the number of indecomposable
1342-avoiding permutations of length $n$ equals that of labeled plane
trees of a certain type on $n$ vertices recently enumerated by Cori,
Jacquard and Schaeffer, which is in turn known to be equal to 
the number of rooted bicubic maps enumerated by Tutte in 1963. Moreover,
$H(x)$ turns out to be algebraic, proving the first nonmonotonic,
longer-than-three instance of a conjecture of Zeilberger and Noonan.
We also prove  that $\sqrt[n]{S_n(1342)}$
converges to 8, so in particular,  $lim_{n\rightarrow
\infty}(S_n(1342)/S_n(1234))=0$.  
\end{abstract}

\maketitle

\section{Introduction} \subsection{Our main results}
In this paper we are going to prove an exact formula for  the number
$S_n(1342)$ of 1342-avoiding permutations of length $n$  showing that  
\[S_n(1342)=\frac{(7n^2-3n-2)}{2}\cdot (-1)^{n-1}
+3\sum_{i=2}^n 2^{i+1}\cdot \frac{(2i-4)!}{i!(i-2)!}\cdot {{n-i+2\choose 2}}
\cdot (-1)^{n-i}, \] by first 
proving that the ordinary generating function $H(x)$ for these numbers
$S_n(1342)$ 
has the following simple form:
\[H(x)=\frac{32x}{-8x^2+12x+1-(1-8x)^{3/2}}.\]

This is the first result we know of which provides an exact formula for
the number $S_n(q)$ of permutations of length $n$ avoiding a given
pattern $q$ that is longer than three and is not 1234. Results
concerning the  case 
of length three can be traced back to two centuries; \cite{catalan}
already makes references to earlier work. The formula for $q=1234$ is
given in \cite{ira}.
Until recently  it has not even been known that
$S_n(1342)<c^n$ for some constant $c$. In \cite{minta} this upper bound
with $c=9$ was proved.  This paper's result pushes down
this $c$ to 8, and proves that it is optimal.  
(Definitions and more background information can be found in the next
subsection). 
 
In our proof, we  are going to present a new link between the enumeration
of permutations avoiding the pattern 1342 and the that of
$\beta(0,1)$-trees, 
a class of labeled trees recently  introduced
 in \cite{cori}.  We will  show
that the  number $I_n(1342)$ of {\em indecomposable}  1342-avoiding
permutations of length $n$ is equal to the number of 
\btrees on $n$ nodes. The set $\deb$ of \btrees 
on $n$ nodes is known \cite{cori} to be  equinumerous to the set of
rooted bicubic maps on $2(n+1)$ vertices,
and an exact formula for the number $t_n$ of these is provided in
\cite{tutte}. Therefore, $I_n(1342)=|\deb|=t_n=3 
\cdot 2^{n-1} \cdot \frac{(2n)!}{(n+2)!n!}$. 
So combinatorially, the number of all 
$n$-permutations avoiding 1342 will be shown to be equal to that of
plane forests on $n$ vertices in which each component is a
$\beta(0,1)$-tree.   To our best knowledge, this is the first time when
permutations avoiding  a given pattern are shown to have such a close
 connection with  some  planar maps, though
recently 2-stack-sortable permutations have been shown to be
equinumerous to nonseparable planar maps \cite{dulucq}
\cite{goulden}.
 
Examining the generating function $H(x)$ we will  prove and
disprove 
several conjectures for the pattern 1342. $H(x)$ turns out to be
algebraic, proving a conjecture of Zeilberger and Noonan \cite{zeno}
for the first time for a nonmonotonic pattern which is longer than
three. We will see that $\sqrt[n]S_n(1342)\rightarrow 8$,
which disproves  a conjecture of Stanley and implies the surprising fact
that $lim_{n\rightarrow \infty}(S_n(1342)/S_n(1234))=0$. 

\subsection{Definitions and Background}
In what follows permutations of length $n$ will be called
$n$-permutations. A  permutation is called {\em indecomposable} if it
 cannot be cut in two parts  so that everything before the cut is
{\em larger} than everything after the cut. Let $q=(q_1, q_2, \ldots,
q_k)$ be a $k$-permutation 
and let $p=(p_1, p_2, \ldots, p_n)$ be an $n$-permutation. We say that
$p$ avoids $q$ if $1\leq i_{q_1}< i_{q_2} <
\ldots < i_{q_k} \leq n$ such that $p(i_1)<p(i_2)<\ldots < p(i_k)$. For
example, $p$ avoids 132 if it cannot be written as $\cdots,a,\cdots,
b,\cdots c,\cdots $ so that $a<c<b$. For another example, a permutation
is 1234-avoiding if it doesn't contain an increasing subsequence of
length 4.

It is a natural and easy-looking question to ask how many
$n$-permutations avoid a given pattern $q$. Throughout this paper, this
number will be denoted by $S_n(q)$.

However, this problem turns out to be very hard. Exact answers have only
been known for the easy case of patterns of length 3 \cite{simion}.
In that case $S_n(q)=C_n$, the $n$-th Catalan number, for any such
pattern $q$.  
If $q$ is longer than three, then the most exact result is due to Regev
\cite{Regev} and says that  for all $n$, $S_n(1234\cdots k)$ 
asymptotically equals $c \frac{(k-1)^{2n}}{n^{(k^2-2k)/2}}$,  where
$c$ is  a constant given by a  multiple integral. 
The major problem of this area is to
prove the conjecture of Wilf and Stanley \cite{wilst} from 1990, stating
that for each pattern $q$ there is an absolute constant $\lambda$ so that
$S_n(q)<\lambda^n$ holds. This has been proven for the case of length 4 and
some longer patterns in \cite{minta}. The general case is, however,
still open. 

Stanley \cite{st3} conjectured that $\sqrt[n]{S_n(q)}$
converges to $(k-1)^2$ where $k$ is the length of $q$.  The results of
this paper will indicate the {\em disproof} of this conjecture for the
pattern 1342, in fact, our formula for $S_n(1342)$ clearly implies that
$\sqrt[n]{S_n(1342)}\rightarrow 8$. This  shows in particular that
$(S_n(1342)/S_n(1234))\rightarrow 0$ when $n\rightarrow \infty$ as 
$\sqrt[n]{S_n(1234)}\rightarrow 9$  by the above results of Regev
\cite{Regev} and Gessel 
\cite{ira}  for monotonic patterns.  We would like to point out the
surprising nature of this discovery: while $S_n(q)={{2n\choose
n}}/(n+1)$ for {\em any} patterns $q$ of length three, for the case of
length four there are sequences $S_n(q)$ that are not only different
from each 
other, but their quotient also converges to 0.

We note that while there are 24 permutation patterns of length 4, for
many of them the sequences $S_n(q)$ are identical. In fact, there are
only three different classes of patterns from this point of view,
\cite{west}, \cite{stankova}, 1342, 1234 and 1324 are representants of
them. 

Recently, attention has been paid to the problem of counting the number
of permutations of length $n$ containing a {\em given number $r$} (as
opposed to 0) of subsequences of a certain type $q$. 
 
The major problem of this field is to describe this function for {\em any
given $r$}, not just for $r=0$. In \cite{zeno} 
Noonan and Zeilberger conjectured that for {\em any} given subsequence
$q$ and for {\em any} given $r$,  the number
of $n$-permutations containing exactly $r$ subsequences of type $q$ is a
\prek function of $n$. Present author has proved this conjecture for any
$r$ when $q=132$. Beyond the case of length 3, however, there have been
no nonmonotic instances solved prior to this paper, not even for the
case of $r=0$. Even in the case of 
length 3,  exact formulae have only been given for $r=1$ (\cite{noonan}
for $q=123$ and \cite{egyes} for $132$). The case of monotonic patterns
and $r=0$ has been handled in \cite{zeilberg}.
Figure 1 shows the present state of research on permutations avoiding
given patterns of length 4, including the contributions of this paper. 

\[ \begin{picture}(200,200)(20,20) 

\put(10,200){\line(1,0){180}} \put(10,80){\line(1,0){180}}
\put(10,20){\line(1,0){180}} \put(10,140){\line(1,0){180}}
\put(10,20){\line(0,1){180}}  \put(130,20){\line(0,1){180}}
\put(70,20){\line(0,1){180}}  \put(190,20){\line(0,1){180}}

\put(25,182){known}\put(32,165){\cite{Regev}}
\put(25,122){known}\put(30,105){\cite{minta}}
\put(25,62){known}\put(30,45){\cite{minta}}

\put(19,205){$S_n(q)<c^n$} \put(78,205){formula} \put(139,205){\prek}
\put(85,182){known}\put(92,165){\cite{ira}} 
\put(85,122){known}\put(76,105){this paper}
\put(85,62){open}
\put(140,182){known}\put(149,165){\cite{zeilberg}}
\put(140,122){known}\put(136,105){this paper}
\put(140,62){open}
\put(-20,45){1324} \put(-20,105){1342} \put(-20,165){1234} 
 \put(78,6){Figure 1}
\end{picture} \]
\vspace*{12pt}

In the following two paragraphs we give a very
brief summary of \prek and algebraic generating functions. The reader
familiar with them can skip these paragraphs.

A function $f:\Nnn \rightarrow \Cnn $ is called $P$-recursive
if there exist polynomials $P_0,P_1,\cdots, P_k\in Q[n]$, with $P_k\neq
0$ so that \begin{equation} \label{rekdef}
P_k(n)f(n)+P_{k-1}(n)f(n+k-1)+\cdots +P_0(n)f(n)=0
\end{equation} for all 
natural numbers $n$. Here \prek stands for ``polynomially
recursive''. The continuous analogue 
of this notion is $d$-finiteness, which stands for ``differentiably
finite''. Let $u(x)\in \Cnn [[x]] $ be a power series. If there
exist polynomials  $p_0(x),p_1(x),\cdots
p_d(x)$  so that $p_d\neq 0$ and
\begin{equation} \label{findef} p_d(x)u^{(d)}(x)+p_{d-1}(x)u^{(d-1)}(x)+\cdots
+p_1(x)u'(x)+p_0(x)u(x)=0, \end{equation} 
then we say that $u$ is $d$-finite. 
(Here $u^{(j)}=\frac{d^ju}{dx^j}$). It is well-known \cite{euro} that  a
function $f(n)$ is \prek if and only if its ordinary  
generating function $u(x)=F(x)=\sum_{n\geq 0}f(n)x^n$ is $d$-finite.
 
Another, smaller class of formal power series is that of {\em algebraic}
series. We say that the series $v(x)\in \Cnn [[x]]$ is {\em algebraic}
if there exist polynomials $p_0(n),p_1(n),\cdots
p_{d-1}(n)$  so that $p_{d-1}\neq 0$ and 

\begin{equation}  v^{d}(x)+p_{d-1}(x)v^{d-1}(x)+\cdots
+p_1(x)v(x)+p_0(x)=0. \end{equation}
 Any algebraic power series is
necessarily $d$-finite. 

As we have mentioned, we are going to prove our results by showing
connections between these permutations and rooted bicubic maps. These
are planar maps with 2-colorable vertices which in addition all have
degree three and a 
distinguished ``root'' edge and face. We will make good use of a class
of labeled rooted trees called \btrees introduced in \cite{cori}. Their
definition is at the beginning of the next section. 

\section{The correspondence between trees and permutations}

\begin{definition} \cite{cori} A rooted plane tree with nonnegative
integer labels $l(v)$ on each of its vertices $v$ is called a
\btree if it satisfies the 
following conditions: 
\begin{itemize} 
\item if $v$ is a leaf, then $l(v)=0$,
\item if $v$ is the root and $v_1,v_2,\cdots, v_k$ are its children,
then $l(v)=\sum_{i=1}^kl(v_k)$,
\item if $v$ is an internal node and $v_1,v_2,\cdots, v_k$ are its children,
then $l(v)\leq 1+ \sum_{i=1}^kl(v_k)$.
\end{itemize}
\end{definition}
A {\em branch} of a rooted tree is a tree whose top is one of the root's
children. Some rooted trees may have only one branch, which doesn't
necessarily mean they consist of a single path.

 We start by treating two special types of \btrees on $n$ vertices.
These cases are fairly simple- they will correspond to 231-avoiding
(resp. 132-avoiding) permutations, but they will be our tools in dealing
with the general case.

First we set up a  bijection $f$ from the set of all
 1342-avoiding $n$-permutations starting with the entry
1 and the set of \btrees on $n$ vertices 
consisting of one single path. In other words, the former is the set of
231-avoiding permutations 
of the set $\{2,3,4\cdots, n\}$. 

So let $p=(p_1p_2\cdots p_n)$ be an 1342-avoiding $n$-permutation 
 so that $p_1=1$. Take an unlabeled tree on $n$ nodes consisting of a
single path 
 and give the label $l(i)$ to its $i$th  node $(1\leq i \leq n-1)$ by
 the following rule: 

$l(i)=\#\{j\leq i $ so that $p_j>p_s$ for at least one $s>i \: \}$.

Finally, let $l(n)=l(n-1)$.
In words, $l(i)$ is the number of entries on the left of $p_i$ (inclusive)
which are larger than at least one entry on the right of $p_i$. 
We note that this way we could define $f$ on the set of {\em all}
$n$-permutations starting with the entry 1, but in that case, as we will
see, $f$ would not be a bijection. (For example, the images of 1342 and
1432 would be identical).

\begin{example} {\em If $p=14325$, then the labels of the nodes of $f(p)$
are, from the leaf to the root,  0,1,2,0,0. } \end{example} 
\begin{lemma} \label{egyenes} 
$f$ is a bijection from the  set of $\deb$ of all \btrees
on $n$ vertices consisting of one single path 
to the set  of 1342-avoiding $n$-permutations starting with
the entry 1. \end{lemma}
\begin{proof}
It is easy to see that $f$ indeed maps into the set of \btrees :
$l(i+1)\leq l(i)+1$ for all $i$ because there can be at most one 
entry counted by $l(i+1)$ and not counted by $l(i)$, namely the entry
$p_{i+1}$. All labels are certainly nonnegative and $l(1)=0$. 

To prove that $f$ is a bijection, it suffices to show that it has an
inverse, that is, for any \btree $T$ consisting of a single path, we can
find 
the only permutation $p$ such that $f(p)=T$. We claim that given $T$, we
can recover the entry $n$ of the preimage $p$. First note that $p$ was
1342-avoiding and started by 1, so any entry on the left of $n$ must
have been smaller than any entry on the right of $n$. In particular, the
node preceding $n$ must have label 0. Moreover, $n$ is the
leftmost entry $p_i$ of $p$ so that  $p_j>0$ for all $j\geq i$ if there
is such an entry at all, and $n=p_n$ if there is none. That is, $n$
corresponds to the node which starts the uninterrupted sequence of
strictly positive labels which ends in the last node, if there is such a
sequence, and corresponds to the last node otherwise. To see this, note
that $n$ is the largest of all entries, so in particular it is always
larger than at least one entry in any nonempty set of entries. 

Once we located where $n$ was in $p$, we can simply delete the node
corresponding to it from $T$ and decrement all labels after it by 1. (If
this means deleting the last node, we just change $l(n-1)$ 
 $l(n-2)$ to satisfy the root-condition). We can indeed do this
because the node preceding $n$ had label 0, and the node after $n$ had a
positive label, by our algorithm to locate $n$. Then we can  proceed
recursively, by finding the position of 
the entries  $n-1,n-2,\cdots,1$ in $p$. This clearly defines the
inverse of $f$, so we have proved that $f$ is a bijection.   
\end{proof}

\begin{lemma} \label{nullak} The number \btrees with all labels equal
to zero is $C_{n-1}$. \end{lemma}
\begin{proof} These \btrees are in fact unlabeled plane trees. We prove
that they are 
in one-to-one correspondence with the 132-avoiding permutations whose
last entry is $n$. Suppose we already know this for all positive
integers $k<n$. Let $T$ be a \btree on $n$ vertices with all labels
equal to 0 and root $r$. Let $r$ have $t$ children, which are, from
left-to-right, at the top of such unlabeled trees $T_1,T_2,\cdots, T_t$
on $n_1,n_2,\cdots, n_t$ 
nodes. Then by induction, each of the $T_i$ corresponds to a
132-avoiding $n_i$-permutation ending with $n_i$. Now add
$\sum_{j=i+1}^t n_j$ to all entries of the permutation $p_i$ associated
with $T_i$, then concatenate all these strings and add $n$ to the end to
get the permutation $p$ associated with $T$. 
This is clearly a bijection as the blocks of the first $n-1$ elements
determine the branches of $T$. \end{proof} 

\begin{example} {\em The permutation 341256 corresponds to the \btree with
all labels equal to 0  shown in  Figure 2.}\normalfont 
\[\begin{picture}(130,130)(30,30)
\put(27,54){\line(1,2){33}} \put(93,54){\line(-1,2){33}}
\put(27,54){\circle*{5}} \put(31,57){3} 
\put(41,81){\circle*{5}} \put(45,84){4}
\put(79,81){\circle*{5}} \put(81,84){2}
\put(93,54){\circle*{5}} \put(95,57){1}
\put(60,118){\line(0,1){17}}
\put(60,118){\circle*{5}} \put(62,116){5}
\put(60,138){\circle*{5}} \put(62,138){6}
\put(32,24){Figure 2}
 \end{picture}\]
\vspace*{12pt}
\end{example}

An easy way to read off the corresponding permutation once we have
its entries written to the corresponding nodes is the well-known  {\em
preorder} 
reading: for every node, first write down the entries associated with its
children from left to right, then the entry associated with the node
itself, and do this recursively for all the children of the node.

Note that such a \btree has only one branch if and only if the
next-to-last element of the indecomposable 132-avoiding $n$-permutation
corresponding to it is $n-1$.

Let's introduce some more notions before we attack the general case. 
An entry of a permutation which is smaller than all the
entries by which it is preceded called a {\em left-to-right minimum}.

\begin{definition} Two $n-$permutations $x$ and $y$ are said to be
 in the same
{\em  class}  if  the left-to-right minima of $ x$
are the same as those of $ y$ and they are in the same positions.
\end{definition} 
\begin{example} {\em 34125 and 35124 are in the same  class since their \ler
minima are 3 and 1, and they are located at the same positions. 3142
and 3412 are not in the same class. } \end{example}

\begin{proposition} Each nonempty class $C$ of $n$-permutations contains
exactly one 132-avoiding permutation. \end{proposition}
\begin{proof} Take all entries which are not \ler minima and fill all
slots between the \ler minima with them as follows: in each step 
place the smallest element which has not been 
placed yet which is larger then the previous left-to-right minimum. The
permutation obtained this way will be clearly 132-avoiding, and it will
be the only one in this class because any time we  deviate
from this procedure, we create a 132-pattern. \end{proof}

\begin{definition} The normalization N(p) of an $n$-permutation $p$ is the
only 132-avoiding permutation in the class $C$ containing $p$.
\end{definition}  
\begin{example} {\em  If $p=32514$, then $N(p)=32415$. } \end{example}
\begin{definition} The normalization $N(T)$ of a \btree $T$ is the
\btree which is isomorphic to $T$ as a plane tree, with all labels equal
to zero. \end{definition}

\begin{proposition} \label{indek} A permutation  $p$ is indecomposable
if and only $N(p)$ is indecomposable. 
\end{proposition}   
\begin{proof} Let $C$ be the class containing $p$, given by the set and
position of its \ler minima.
It is clear that if $p\in C$ is decomposable, then the only way to cut it in
two parts (so that everything before the cut is larger than everything
after the cut) is to cut it immediately before a \ler minimum $a$. Now if
there is a \ler minimum $a$  so that it is in the $n+1-a$th
position, then all entries which are larger than $a$ must be placed on
the left of $a$ and so all such permutations in $C$ are decomposable. If there
is no such $a$, then for all \ler minima $m$, there will be an entry $b$
so that 
$m<b$ and $b$ is on the right of $m$ and so  permutations in $C$ will not
be decomposable. \end{proof}

\begin{corollary} \label{vegzodes} If $p$ is an indecomposable
$n$-permutation, then 
$N(p)$ always ends with the entry $n$. \end{corollary}
\begin{proof} Note that the only way for a 132-avoiding $n$-permutation to
be indecomposable is for it  to end with $n$. Then the statement follows from
Proposition \ref{indek}. \end{proof} 

Now we are in a position to prove our theorem about the number of
indecomposable 1342-avoiding permutations. 

\begin{theorem}\label{bijek} The number of indecomposable
$n$-permutations which 
avoid the pattern 1342 is \[I_n(1342)=t_n=3 \cdot 2^{n-1} \cdot
\frac{(2n)!}{(n+2)!n!} \]
\end{theorem} 
\begin{proof} We are going to set up a bijection $F$ between these
permutations and $\deb$. This will be an extension of the bijection $f$
of lemma \ref{egyenes}. As the size of $\deb$ is known to be equal to
$t_n$ \cite{cori}, this will prove our claim.

Let $p$ be and indecomposable 1342-avoiding $n$-permutation. Take
$N(p)$. By corollary \ref{vegzodes} its last element is $n$. Define
$F(N(p))$ to be the \btree $S$ associated to $N(p)$ by the bijection of
lemma \ref{nullak}. Now write the entries of $p$ to the nodes of $S$ so
that for all $i$, the $p_i$ is written to the node where $N(p)_i$ was
written in $S$. In particular, the \ler minima remain unchanged. 
Figure 3 shows how we associate the entries of the permutation $361542$
to the nodes of $N(T)$, which is the image of $N(p)=341256$.

\[\begin{picture}(130,130)(30,30)
\put(27,54){\line(1,2){33}} \put(93,54){\line(-1,2){33}}
\put(27,54){\circle*{5}} \put(31,57){3} 
\put(41,81){\circle*{5}} \put(45,84){6}
\put(79,81){\circle*{5}} \put(81,84){5}
\put(93,54){\circle*{5}} \put(95,57){1}
\put(60,118){\line(0,1){17}}
\put(60,118){\circle*{5}} \put(62,116){4}
\put(60,138){\circle*{5}} \put(62,138){2}
\put(30,23){Figure 3}
 \end{picture}\]
\vspace*{12pt}

Now we are going to define the label of each node for this new \btree $T$
and  obtain $F(p)$  this way. (As an unlabeled tree, $T$ will be
isomorphic to $S$, but its labels will be different). Denote $i$ the
$i$th node of  $T$ in the preorder reading, thus $p_i$ is the $i$th
entry of $p$, (which is therefore associated to node $i$), while $l(i)$
is the label of this  node. We say that $p_i$ {\em beats} $p_j$ if there
is an element $p_h$ so that $p_h,p_i,p_j$ are written in this order and
they form a 
132-pattern. Moreover, we say that $p_i$ {\em reaches} $p_k$ if there
is a subsequence $p_i,p_{i+a_1},\cdots p_{i+a_t},p_k$ of entries so that
$i<i+a_1<i+a_2<\cdots <i+a_t<k$ and that any entry in this subsequence
beats the next one. For example, in the permutation 361542, the entry
6 beats 5 and 4, 5 beats 4 and 2, and 4 beats 2, while 6 reaches 2 (of
course, each entry reaches all those elements it beats, too).
 Then let

$ l(i)= \# \{j$ descendants of $i$ (including $i$ itself) so that
there is at least one  $k>i$ for which $p_j$ reaches $p_k \},$

and let $F(p)$ be the \btree defined by these labels. (Recall that a
descendant of $i$ is an element of the tree whose top 
element is $i$). First, it is easy to see that $F$ indeed maps into the
set of $\beta(0,1)$-trees: if $v$ is an internal node and
$v_1,v_2,\cdots, v_k$ are its children,  
then $l(v) \leq 1+ \sum_{i=1}^k l(v_k)$  because there can be at most one 
entry counted by $l(v)$ and not counted by any of its children's label,
namely $v$ itself. All labels are certainly nonnegative and all leaves,
that is, the \ler minima, have label 0.

If $p=361542$, then $F(p)$ is the \btree shown in Figure 4b. 

\[\begin{picture}(230,130)(30,30)

\put(27,54){\line(1,2){33}} \put(93,54){\line(-1,2){33}}
\put(27,54){\circle*{5}} \put(31,57){3} 
\put(41,81){\circle*{5}} \put(45,84){6}
\put(79,81){\circle*{5}} \put(81,84){5}
\put(93,54){\circle*{5}} \put(95,57){1}
\put(60,118){\line(0,1){17}}
\put(60,118){\circle*{5}} \put(62,116){4}
\put(60,138){\circle*{5}} \put(62,138){2}
\put(30,20){Figure 4a}

\put(127,54){\line(1,2){33}} \put(193,54){\line(-1,2){33}}
\put(127,54){\circle*{5}} \put(131,57){0} 
\put(141,81){\circle*{5}} \put(145,84){1}
\put(179,81){\circle*{5}} \put(181,84){1}
\put(193,54){\circle*{5}} \put(195,57){0}
\put(160,118){\line(0,1){17}}
\put(160,118){\circle*{5}} \put(162,116){3}
\put(160,138){\circle*{5}} \put(162,138){3}
\put(130,20){Figure 4b}
 \end{picture}\] 
\vspace*{12pt}

To prove that $F$ is a bijection, it suffices to show that it has an
inverse, that is, for any \btree $T\in \deb$, we can find
the only permutation $p$ so that $F(p)=T$. We again claim that given
$T$, we can recover the node $j$ which has the entry $n$ of the preimage $p$
associated to it, and so we can recover  the position of $n$ in the preimage.

\begin{proposition} \label{pozitiv} Suppose $p_n\neq n$, that is, $n$ is
not associated to the root vertex. Then each
ancestor of $n$, including $n$ itself, has a positive label. If
$p_n=n$, then $l(n)=0$ and thus there is no vertex with the above property.
\end{proposition}
\begin{proof}
If $p_n=n$, then nothing beats it, thus $p_n=0$. Suppose $p_n$ is not
the root vertex. 

To prove our claim  it is enough to show that for any node $i$ which is an
an ancestor of $j$, there is an entry $p_k$ so that $p_k$ is an ancestor
 of
$p_i$ and $n=p_j$ reaches $k$. Indeed, this would imply that the entry $p_j=n$
is counted by the label $l(i)$ of $i$. Now let $a_m=p_1>a_2>\cdots
a_1=1$ be the \ler minima of $p$ so that $n$ is located between $a_r$
and $a_{r+1}$. Then $n$ certainly beats all elements located between $a_r$ and
$a_{r+1}$ as $a_r$ can play the role of 1 in the 132-pattern. Clearly,
$n$ must beat at least one entry $y_1$ on the right of $a_{r+1}$ as well,
otherwise $p$ would be decomposable by cutting it right before
$a_{r+1}$.  If $y_1$ is on the right of $i$, then we are done. If not,
then  $y_1$ must beat at least one entry $y_2$ which is on the
other side of $a_{r_1+1}$, where $y$ is located between $a_{r_1}$ and
$a_{r_1+1}$ for the same reason, and so on. This way we get a subsequence
$y_1, y_2, \cdots $ so that $n$ reaches each of the $y_t$, and this
subsequence eventually gets to the right of $i$, since in each step we
bypass at least one \ler minimum. Thus the proposition is proved.
\end{proof}

\begin{proposition} Suppose $p_n\neq n$. Then $n$ is the {\em leftmost}
entry of $p$ which has the property that each of its ancestors has a
 positive label. \end{proposition}

\begin{proof} Suppose $p_k$ and $n$ both have this property and that
$p_k$ is on the left of $n$. (If there are several candidates for the
role of $p_k$, choose the rightmost one). If $p_k$ beats an element $y$
on the right 
of $n$ by participating in the 132-pattern $x\:p_k\:y$, then
$x\:p_k\:n\:y$ is a 1342-pattern, which is a contradiction. So $p_k$
does not beat such an element $y$. In other words, all elements after
$n$ are smaller than all elements before $p_k$. Still, $p_k$ must reach
elements on the right of $n$, thus it beats some element $v$ between
$p_k$ and $n$. This element $v$ in turn beats some element $w$ on the
right of $n$ by participating in some 132-pattern $t\:v\:w$. However,
this would imply that $t\:v\:n\:w$ is a 1342-pattern, a contradiction,
which proves our claim. 
\end{proof}

Therefore, we can recover the entry $n$ of $p$ from $T$. Then we can
 proceed as in the proof of Lemma 1, that is, 
 just delete $n$, subtract 1 from the labels of its ancestors and
 iterate this procedure to get $p$. If any time during this procedure
we find that the current root is associated to the maximal entry that
hasn't been associated to other vertices yet, and the tree has more than
one branch at this moment, then deleting the root vertex will split the
tree into smaller trees. Then we continue the same procedure on each of
them. The set of the entries associated to each of these smaller trees is
uniquely determined because $T$ as an unlabeled tree determines the \ler
minima of $p$.  Therefore, we can always recover $p$ in this way from
$T$. This proves that $F$ is a bijection.
  
 Thus we have set up a bijection between the set
of indecomposable 1342-avoiding $n$-permutations and $\deb$. We know
from \cite{cori} that $|\deb|= t_n=3 
\cdot 2^{n-1} \cdot \frac{(2n)!}{(n+2)!n!}$ and therefore the Theorem is
proved.  
\end{proof} 

Note that in particular, $F$ maps 132-avoiding permutations into \btrees
with all labels equal to 0 and permutations starting with the entry 1
into \btrees consisting of a single path.

\begin{corollary} $S_n(1342)$ equals the number of plane forests on $n$
vertices in which
each component is a $\beta(0,1)$-tree. \end{corollary}

\section{Enumerative results} 
Tutte \cite{tutte} has obtained the numbers $t_n$ by first computing a
translate of their generating function
\begin{equation} \label{genfugg} F(x)=\sum_{n=1}^{\infty} 3 \cdot
2^{n-1}  \cdot
\frac{(2n)!}{(n+2)!n!}x^n= \frac{8x^2+12x-1+(1-8x)^{3/2}}{32x}.
\end{equation} 
By  theorem \ref{bijek}, the coefficients of this generating function are the
numbers $I_n(1342)$.  Therefore, the generating function of {\em all}
1342-avoiding permutations is given by the following theorem.
\begin{theorem} Let $s_n=S_n(1342)$ and let 
$H(x)=\sum_{n=0}^{\infty}s_nx^n$. Then
\begin{equation} \label{altal} 
H(x)=\sum_{i\geq 0}F^i(x)=\frac{1}{1-F(x)}=
\frac{32x}{-8x^2+12x+1-(1-8x)^{3/2}}. \end{equation}
\end{theorem}
\begin{proof}  Any 1342-avoiding permutation 
has a unique decomposition into indecomposable permutations. This can
consist of $1,2,\cdots$ blocks, 
implying that $s_n=\sum_{i=1}^n t_is_{n-i}$, and the statement follows. 
\end{proof}

\begin{theorem}\label{explicit} For all $n\geq 0$ we have 
\begin{equation} 
\def\eqalign#1{\null\,\vbox{\openup\jot
  \ialign{\strut\hfil$\displaystyle{##}$&$\displaystyle{{}##}$\hfil
      \crcr#1\crcr}}\,}
\label{eformula}
\medmuskip1mu
\hskip-1em
\eqalign{S_n(1342)&=S_n(1342)\cr&=\frac{(7n^2-3n-2)}{2}\cdot (-1)^{n-1} 
+3\sum_{i=2}^n 2^{i+1}\cdot \frac{(2i-4)!}{i!(i-2)!}{{n-i+2\choose 2}}
(-1)^{n-i}.}
\hskip-1em
\end{equation} 
\end{theorem}
\begin{proof} Multiply both the numerator and the denominator of $H(x)$
by $(-8x^2+20x+1)+(1-8x)^{3/2}$. After simplifying we get
\begin{equation} \label{egyszeru}
H(x)=\frac{(1-8x)^{3/2}-8x^2+20x+1}{2(x+1)^3}. \end{equation}
As $(1-8x)^{3/2}=1-12x+\sum_{n\geq 2} 3\cdot 2^{n+2}x^n
\frac{(2n-4)!}{n!(n-2)!}$, formula (\ref{egyszeru}) implies our claim.
\end{proof}

So the first few values of $S_n(1342)$ are 
1,2,6,23,103,512,2740,15485,91245,555662. 
In particular, one sees easily that the expression on the right hand
side of (\ref{eformula}) is dominated by the last summand; in fact, the
alternation in sign assures that this last summand is larger than the
whole right hand side if $n\geq 8$. As
$\frac{(2n-4)!}{n!(n-2)!}<\frac{8^{n-2}}{n^{2.5}}$ by Stirling's formula, we
have proved the following exponential upper bound for $S_n(1342)$.
\begin{corollary}\label{upper} For all $n$, we have
$S_n(1342)<8^n$. \end{corollary}

It is straightforward to check that the numbers $I_n=t_n$ satisfy the
following recurrence 
\begin{equation}\label{rekurziv} t_n=(8n-4)t_{n-1}/(n+2). \end{equation}

In particular, $\sqrt[n]{t_n}\rightarrow 8$. 
Using this formula we can disprove a conjecture of Stanley claiming that
for all permutation patterns $q$ of length $k$, the sequence
$\sqrt[n]{S_n(q)}$ converges to $(k-1)^2 $.  

\begin{theorem} $\sqrt[n]{s_n}=\sqrt[n]{S_n(1342)}\rightarrow 8$ when
$n\rightarrow \infty $. \end{theorem}

\begin{proof} This is true as clearly $t_n\leq s_n < 8^n$  by Corollary
\ref{upper}  and we know from (\ref{rekurziv}) that 
$\sqrt[n]{t_n}\rightarrow 8$ if $n\rightarrow \infty$. 
\end{proof}

\begin{corollary} 
 $lim_{n\rightarrow \infty}(S_n(1342)/S_n(1234))=0$. \end{corollary}
\begin{proof} Follows from $lim_{n\rightarrow \infty}
\sqrt[n]{S_n(1234)}=9$ \cite{ira} \cite{Regev}. \end{proof}

This Corollary certainly implies that $S_n(1342)<S_n(1234)$ if $n$ is
large enough. However, using the formulae of Theorem \ref{explicit} and
\cite{ira}, we can easily show that this is true for {\em all} $n\geq
6$. (This has recently been shown by a long argument in
\cite{minta}).

\begin{corollary} For all $n\geq 6$, we have
$S_n(1342)<S_n(1234)$. \end{corollary} 
\begin{proof} It is known \cite{ira} that \begin{equation} \label{irag}
S_n(1234)=2\cdot \sum_{i=0}^n {{2i\choose i}}\cdot {{n\choose i}}^2
\cdot \frac{3k^2+2k+1-n-2kn}{(k+1)^2(k+2)(n-k+1)} \end{equation}
One sees easily that the dominant summand is the one with $i=2n/3$, and
that this summand is much larger than the last (and dominant) summand in
(\ref{eformula}) if $n\geq 9$. The proof then follows by checking the
values of $S_n(1342)$ and $S_n(1234)$ for $n\leq 8$. \end{proof}

Formula (\ref{altal}) enables us to prove that the sequence $S_n(1342)$
is \prek in $n$, solving an instance of the conjecture of Zeilberger and
Noonan mentioned in the Introduction. Indeed, $H(x)$ is certainly algebraic,
thus in particular, it is $d$-finite and therefore $S_n(1342)$ is \prek
as claimed. So we have proved the following theorem.
 
\begin{theorem} The sequence $S_n(1342)$ is \prek in $n$. Furthermore,
its generating function
$H(x)$ is algebraic and its only irrationality is $\sqrt{1-8x}$.
\end{theorem}  

\section*{Acknowledgement}
I am grateful to my research advisor Richard Stanley and to Robert Cori
who indicated me very useful references, as well as to Gilles Schaeffer
who sent me the preprint \cite{cori}. I am also indebted to Sergey Fomin
for his helpful remarks and suggestions.

\end{document}